----------
X-Sun-Data-Type: default
X-Sun-Data-Description: default
X-Sun-Data-Name: grodr.tex
X-Sun-Charset: us-ascii
X-Sun-Content-Lines: 766

\input amstex
\input amsppt.sty
\magnification\magstep1

\def\ni\noindent
\def\sbs{\subset}

\def\R{\text{\bf R}}

\def\H{\text{\bf H}}
\def\N{\text{\bf N}}

\def\e{\epsilon}

\def\sU{\Cal U}
\def\sV{\Cal V}

\hoffset= 0.0in
\voffset= 0.0in
\hsize=32pc
\vsize=43pc
\baselineskip=24pt
\NoBlackBoxes
\topmatter
\author
A.N. Dranishnikov
\endauthor

\title
On large scale properties of manifolds
\endtitle
\date{10.25.99}
\enddate
\abstract
We show that a space with a finite asymptotic dimension is embeddable in
a non-positively curved manifold. Then 
we prove that if a uniformly contractible manifold $X$ is uniformly embeddable in $\R^n$ or non-positively curved $n$-dimensional simply connected manifold
then $X\times\R^n$ is integrally hyperspherical. 
If a 
uniformly contractible manifold $X$ of bounded geometry is uniformly
embeddable into a Hilbert space, then $X$ is stably integrally hyperspherical.
\endabstract

\thanks The author was partially supported by NSF grant DMS-9971709.
\endthanks

\address Pennsylvania State University, Department of Mathematics, 
218 McAllister Buil\-ding, University Park, PA 16802 , USA
\endaddress

\subjclass Primary 20H15 
\endsubjclass

\email  dranish\@math.psu.edu
\endemail

\keywords  hyperspherical manifold, uniform embedding, 
scalar curvature, Novikov conjecture.
\endkeywords
\endtopmatter

\document
\head \S1 Introduction \endhead

Gromov introduced several notions of largeness of Riemannian manifolds.
For example, a manifold $X$ which is a universal cover of a closed aspherical manifold $M^n$ with the fundamental group $\Gamma=\pi_1(M^n)$, supplied with a 
$\Gamma$-invariant metric is large in a sense that it is  uniformly contractible. We recall that $X$ is {\it uniformly contractible} if there is a function $S(r)$ such that every ball $B_r(x)$ of the radius $r$ centered at
$x$ is contractible to a point in the ball $B_{S(r)}(x)$ for any $x\in X$.
For the purpose of the Novikov conjecture and other related conjectures, it is important to show that universal covers of aspherical manifolds are also large in some cohomological sense. The weakest such property is called hypersphericity.
\proclaim{Definition [G-L]}
An $n$-dimensional manifold $X$ is called hyperspherical if for every $\epsilon>0$
there is a $\epsilon$-contracting proper map $f_{\epsilon}:X\to S^n$ of nonzero
degree onto the standard unit $n$-sphere.
\endproclaim
Here a continuous map $f:X\to S^n$ is proper if it has only one unbounded preimage. Gromov and Lawson proved the following [G-L] 
\proclaim{Gromov-Lawson's Theorem}
 An aspherical manifold with a
hyperspherical universal cover cannot carry a metric of a positive scalar
curvature. 
\endproclaim
The natural question appeared [G2]:
\proclaim{Problem 1}
Is every uniformly contractible manifold hyperspherical?
\endproclaim
Above we defined a notion of the rational hypersphericity. One can define
integral hypersphericity by taking maps of degree one. Also the $n$-sphere
can be replaced by $\R^n$. In that case we obtain the notion of a  {\it hypereuclidean}
manifold. In the integral case these two notions seems coincide. Also in the integral
case the above question of Gromov has a negative answer [D-F-W]. Still there is
no candidates for a rational counterexample and even new examples of Gromov
[G4] leave a possibility for an affirmative answer. 

The following definition is due to Gromov [G1]:
\proclaim{Definition}
A map $f:X\to Y$ between metric spaces is called {\it large scale uniform embedding} if there are two functions $\rho_1,\rho_2:[0,\infty)\to[0,\infty)$
tending to infinity such that

$\rho_1(d_X(x,x'))\le d_Y(f(x),f(x'))\le \rho_2(d_X(x,x'))$.
for all $x,x'\in X$.
\endproclaim
In Section 2 of this paper we prove the following:
\proclaim{Theorem 1}
Suppose that a uniformly contractible manifold $X$ admits a large scaale
uniform embedding into $\R^n$. Then $X\times\R^n$ is integrally hyperspherical.
\endproclaim
\proclaim{Corollary}
Let $M$ be a closed aspherical manifold and assume that the fundamental group
$\Gamma=\pi_1(M)$ admits a coarsely uniform embedding into $R^n$ as a metric space with the word metric. Then $M$ cannot carry a metric with a positive scalar curvature.
\endproclaim
 Some results of that type were already known to Gromov 
(see his remark (b), page 183 of [G3]). Also the Corollary follows from the theorems
of Yu [Yu1],[Yu2].
\proclaim{The First Yu's Theorem} If a finitely presented group $\Gamma$ 
has a finite asymptotic dimension, in particular if it is coarsely uniformly embeddable in
$\R^n$, then the coarse Baum-Connes conjecture holds for $\Gamma$.
\endproclaim

\proclaim{The Second Yu's Theorem} If a finitely presented group $\Gamma$ can be
uniformly in large scale sense embedded into a Hilbert space, then the coarse Baum-Connes conjecture holds for $\Gamma$.
\endproclaim
We note that The Second Yu's theorem implies The First [H-R2].
Open Riemanian manifolds which obey the coarse Baum-Connes conjecture
are of course large in some refined sense. This largeness is relevant to the 
   hypersphericity
or the hypereuclidianess but it is different. The best we can say [Ro1] that integrally
(rationally) hypereuclidean manifolds satisfy the monicity part of the coarse
(rational) Baum-Connes conjecture. We note that the monicity part of the coarse
Baum-Connes is sufficient for the Gromov-Lawson conjecture about nonexistence
of a positive scalar curvature metric on a closed aspherical Reimannian manifolds.

Both Yu's theorems hold for all proper metric spaces with bounded geometry. We recall that a metric space $X$ is called having a {\it bounded geometry} if for every
$\epsilon>0$ and every $r>0$ there is $c$ such that the $\epsilon$-capacity of
every $r$-ball $B_r(x)$ does not exceed $c$. The latter means that a ball
$B_r(x)$ contains no more than $c$ $\epsilon$-disjoint points.
\proclaim{Definition}
An $n$-dimensional manifold $X$ is called stably (integrally) hyperspherical if for any
$\epsilon>0$ there is $m$ such that $X\times\R^m$ admits an $\epsilon$-contracting map of nonzero degree (of degree one) onto the unit $(n+m)$-sphere $S^{n+m}$.
\endproclaim
We prove the following
\proclaim{Theorem 2}
Suppose that $X$ is a uniformly contractible manifold with bounded geometry
and assume that $X$ admits a coarsely uniform embedding into a Hilbert space.
Then $X$ is stably integrally hyperspherical.
\endproclaim
It is unclear whether the stable hypersphericity implies the Gromov-Lawson
conjecture. A positive answer to the following problem would give a simple
argument for the implication: Stable Hypersphericity $\Rightarrow$ 
Gromov-Lawson.
\proclaim{Problem 2}
Does there exist a positive constant $c>0$ such that the $K$-area of the
unit $2n$-sphere is greater than $c$ for all $n$ ?
\endproclaim
We recall the definition of $K$-area from [G3]. For an even dimensional
Riemannian manifold $M$ the $K$-area is

\centerline{$K-area(M)=(\inf_X\{\|R(X)\|\})^{-1}$}
where the infimum is taken over all bundles $X$ over $M$ with some of the Chern numbers nonzero and with unitary connections, $R(X)$ is the curvature of $X$
equipped with the operator norm, i.e. $\|A\|=\sup_{\|x\|=1}\|Ax-x\|$.

The other way to establish the implication:  Stable Hypersphericity $\Rightarrow$ 
Gromov-Lawson would be extending the notion of $K$-area on loop spaces and 
showing that the $K$-area of $\Omega^{\infty}\Sigma^{\infty}S^n$ is greater than
zero.

One of the objectives of this paper is to give  elementary proofs of
Yu's theorems with replacing the coarse Baum-Connes conjecture by the 
hypersphericity. It is accomplished with the First Yu's Theorem. In Section 3 
we  prove The Embedding Theorem for asymptotic dimension. Then 
a slightly more general version of Theorem 1 completes the proof.
The Theorem 2 can be considered as a version of the second Yu's Theorem.
Neverertheless the same set of corollaries 
 would follow from it in the case of 
positive answer to Problem 2.
Otherwise to obtain  a proper analog of the Second Yu's theorem on this way one should do more elaborate 
differential geometry in the argument in order to avoid stabilizing with $\R^m$.

\head \S2 Proofs of Theorems 1 and 2\endhead

First we note that in both theorems it suffices to consider the case when
$X$ is isometrically embedded in $\R^n$ ( or $l_2$). Since a contractible 
$k$-manifold $X$ is homeomorphic to $\R^{k+1}$ after crossing with $\R$, without 
loss of generality we may assume that $X$ is homeomorphic to $\R^k$. Fix a 
homeomorphism $h:\R^k\to X$. We denote by $S_r^{k-1}$ the standard sphere in
$\R^k$ of radius $r$ with the center at $0$. Note that the family $h(S^{k-1}_r)$
tends to infinity in $X$ as $r$ approaches infinity. Denote by $N_{\lambda}(A)$
the $\lambda$-neighborhood of $A$ in an ambient space $W$.
By $B^m_r$ we denote the standard $r$-ball in $\R^m$.
\proclaim{Lemma 1} Let $X$ be a uniformly contractible manifold with bounded
geometry, $X$ is homeomorphic to $\R^k$ and $X$ is isometrically embedded in
a metric space $W$. Then for any $\lambda>0$ there is $r>0$ such that the
neighborhood $N_{\lambda}(h(S^{k-1}_r))$ in $W$ admits a retraction onto
$h(S^{k-1}_r)$.
\endproclaim
\demo{Proof}
Let $S:\R_+\to\R_+$ be a monotone contractibility function on $X$. 
Clearly, $S(t)\ge t$.
Since $X$ is
a space of bounded geometry, there is a uniformly bounded cover $\sU$ on $X$ of
finite multiplicity $m$ and with the Lebesgue number $>4\lambda$
(see [H-R1] or [Dr]). If we squeeze
every element $U\in\sU$ by $\lambda$ we get a cover $\sU'$ with the Lebesgue
number $3\lambda$ For every $U\in\sU'$ we consider the open $\lambda$-neighborhood $ON_{\lambda}(U)$ in $W$. Then the cover
$\tilde\sU=\{ON_{\lambda}(U)\mid U\in\sU'\}$ has the multiplicity $\le m$.

Let $d$ be an upper bound for the diameters of elements of the cover
$\tilde\sU$. Define $T(t)=2S(4t)$. Let $T^l$ mean $l$ times iteration of $T$.
We take $r$ such that $d(h(0),h(S^{k-1}_r))\ge T^{m+k}(d)$.
Denote by $\bar\sU$ the restriction of $\tilde\sU$ over $h(S^{k-1}_r)$.
Note that $\bar\sU$ covers the neighborhood $N_{\lambda}(h(S^{k-1}_r))=N$.
Let $\nu:N\to N(\bar\sU)$ be a projection to the nerve of the cover $\bar\sU$.
We define a map $\phi:N(\bar\sU)\to X\setminus\{h(0)\}$ such that the restriction $\phi\circ\nu\mid_{h(S^{k-1}_r)}$ is homotopic to the identity map
$id_{h(S^{k-1}_r)}$. Then by the Homotopy Extension Theorem there is an extension $\beta:N\to X\setminus\{h(0)\}$ of the identity map $id_{h(S^{k-1}_r)}$. Let $\gamma:\R^k\setminus\{0\}\to S^{k-1}_r$ be a retraction. We define a retraction $\alpha:N\to h(S^{k-1}_r)$ as $h\circ\gamma\circ h^{-1}\beta$.

We define $\phi$ on the $l$-skeleton of $N(\bar\sU)$ by induction on $l$ in such a way that the diameter of the image $\phi(\sigma^l)$ of every $l$-simplex
$\sigma^l$ does not exceed $T^l(d)$. To do that, we choose points $x_U\in U\cap
h(S^{k-1}_r)$ for every $U\in\bar\sU$ and define $\phi(v_U)=x_U$ for every vertex $v_U$ in $N(\bar\sU)$ corresponding to an open set $U\in\bar\sU$.
For every edge $[v_U,v_{U'}]$ we define $\phi$ on it in such a way that
$diam\phi([v_U,v_{U'}])\le S(d(x_U,x_{U'}))\le S(2d)\le T(d)$. 
Assume that $\phi$ is defined on the $l$-skeleton with the property that 
$diam\phi(\sigma)\le T^l(d)$ for all simplices $\sigma$. 
If the boundary is connected, then for arbitrary
$l+1$-dimensional simplex $\Delta$ the image of the boundary 
$\phi(\partial\Delta)$ has the diameter $\le 4T^l(d)$. Then we 
can extend $\phi$
over $\Delta$ with the diameter $\phi(\Delta)\le 2S(4T^l(d))=T(T^l(d))=T^{l+1}(d)$. A map $\phi$ constructed on this way has the
property that $d(x,\phi\nu(x))\le 2(d+T^m(d))\le T^{m+1}(d)$.

Consider a small (with mesh smaller than $\frac{1}{4}T^{m+1}(d)$\ )
triangulation $\tau$ on $h(S^{k-1}_r)$ and the cellular complex
$h(S^{k-1}_r)\times I$ defined by that structure. Using induction one can extend
the map $id_{h(S^{k-1}_r)\times\{0\}}\coprod\phi\nu_{h(S^{k-1}_r)\times\{1\}}:
h(S^{k-1}_r)\times\{0,1\}\to X\setminus\{h(0)\}$ to a map $H:h(S^{k-1}_r)\times I\to X\setminus\{0\}$ with the diameter of the image of every $i$-dimensional
cell less than $T^{m+i+1}(d)$. Then $H(h(S^{k-1}_r))\subset N_{T^{m+k}(d)}(h(S^{k-1}_r))\subset X\setminus\{h(0)\}$  by the choice of $r$.
Thus, $H$ is a required homotopy.\qed
\enddemo
\proclaim{Lemma 2}
Let $M$ be a closed smooth $l$-dimensional submanifold in the euclidean space
$\R^n$ with a trivial tubular $\epsilon$-neighborhood $N_{\epsilon}(M)$.
Then for any number $d>\epsilon$ there is a number $\mu$ such that the diagonal
embedding $j=(1_M\Delta\mu 1_M):M\to\R^n\times\R^n$ has a regular neighborhood
$N$, $\nu:N\to j(M)$ such that:
\roster
\item{} $N$ admits a contracting map $h:N\to B^{2n-l}_d$ which is a homeomorphism
on every fiber $\nu^{-1}(y)$;
\item{} $pr_1(N)\subset N_{2d}(M)$ where $pr_1:\R^n\times\R^n\to\R^n$ is the
projection onto the first factor.
\endroster
\endproclaim
\demo{Proof}
Let $q:N_{\epsilon}(M)\to B_{\epsilon}^{n-l}$ be a trivialization of the 
tubular neighborhood $N_{\epsilon}(M)$. Let $\lambda$ be its Lipschitz constant.
Take $\mu=\frac{\lambda d}{\epsilon}$. Also by $\mu$ we denote the map
$\mu:\R^n\to\R^n$ which is a multiplication of vectors by $\mu$. We extend
the embedding $j:M\to\R^n\times\R^n$ to the map $\bar j:\R^n\to\R^n\times\R^n$
defined as $\bar j(x)=(x,\mu x)$. Thus, the map $\bar j$ is a homothetic transformation
of $\R^n$ to the space $L=\{(x,\mu x)\mid x\in\R^n\}$ with the homothety 
coefficient equal to $\sqrt{1+\mu^2}$. Therefore $j(M)$ admits a trivial tubular 
$\delta$-neighborhood $N'_{\delta}$ in $L$ with $\delta=\epsilon\sqrt{1+\mu^2}$.
We define $N$ as the product $N'_{\delta}\times B^n_d$ isometrically realized in
$L\times L^{\perp}$ where $L^{\perp}$ is the orthogonal complement of $L$ in
$\R^n\times\R^n$. We consider the map 
$h_1=\frac{d}{\epsilon}q\circ\bar j^{-1}
\mid_{N'_{\delta}}:N'_{\delta}\to B^{n-l}_d$. Note that
$\frac{d}{\epsilon}\lambda (1/\sqrt{1+\mu^2})$ is a Lipschitz constant
for $h_1$. Therefore $1=\frac{d}{\epsilon}\frac{\lambda}{\mu}$ is also a
Lipshitz constant for $h_1$. Hence the map $\bar h=h_1\times id_{B^n_d}:
N\to B^{n-l}_d\times B^n_d$ is a short map. Let $p: B^{n-l}_d\times B^n_d\to
B^{2n-d}_d$ be the natural radial projection. We define $h=p\circ\bar h$.
Clearly, the condition 1) holds. Now let $y\in N$ be an arbitrary point, show that
$pr_1(y)\in N_{2d}(M)$. First $y$ can be presented as $j(x)+w_1+w_2$ where
$x\in M$, $w_1\in B^{n-l}_{\delta}\cong \nu_1^{-1}(j(x))$, $\nu_1:N'_{\delta}\to
j(M)$ is the natural projection, and $w_2\in B^n_d\subset L^{\perp}$.
Then $pr_1(j(x)+w_1+w_2)=pr_1(j(x))+pr_1(w_1)+pr_1(w_2)=x+u+pr_1w_2$ where
$w_1=\bar j(u)=(u,\mu u)$ and $u$ is a normal vector to $M$  at point $x$ of the length
$\le\epsilon$. Then $\|pr_1(y)-x\|=\|u+pr_1(w_2)\|\le\|u\|+\|w_2\|\le
\epsilon+d\le 2d$. Hence $pr_1(y)\in N_{2d}(M)$.\qed
\enddemo

\

\demo{Proof of Theorem 1} Let $\dim X=k$.
For every $d>0$ we construct a submanifold $V\subset X\times\R^n$ and a short
map of degree one $f:(V,\partial V)\to (B^{k+n}_d,\partial B^{k+n}_d)$. 
Clearly, this would imply the integral hypersphericity of $X\times\R^n$.
By Lemma 1 for large enough $r$ there is a retraction of the $2d$-neighborhood
$N_{2d}(h(S_r^{k-1}))$ onto a curved $(k-1)$-sphere $h(S_r^{k-1})$. Let
$\alpha$ be the retraction. We may require that $\alpha$ as well as $h$ are
smooth maps.
We assume that $\R^n\subset S^n$
is compactified to the $n$-sphere and assume that $f:S^n\to B^k_r$ is a smooth
extension of $h^{-1}\circ \alpha$. Let $x_0$ be a regular value of $f$ and assume that $f(\infty)\ne x_0$. Then the fiber $M=f^{-1}(x_0)$ is a closed
$(n-k)$-dimensional manifold which admits a
trivial tubular neighborhood $N_{\epsilon}(M)$ for some $\epsilon>0$.
We may assume that $\epsilon<d$.
Then we apply Lemma 2 to obtain an embedding $j:M\to\R^n\times\R^n$ with a
regular neighborhood $N$ with the properties (1)-(2) of the lemma.
In view of condition (2) of Lemma 2
for large enough $R$ the neighborhood $N$ is contained in $N_{2d}(N)\times B^n_R$. Hence the boundary $\partial (h(B^k_r)\times B^n_R)=h(S^{k-1}_r)\times
B^n_R\cup h(B^k_r)\times\partial B^n_R$ does not intersect $N$.
Note that the manifold $M$ is linked in $\R^n$ with $h(S_r^{k-1})$ with the
linking number one. Also $M$ is linked with $\partial (h(B^k_r)\times B^n_R)$
with the linking number one. 
Since $j(M)$ is homotopic to $M$ inside $N_{2d}(M)\times
Int(B^n_R)$, it follows that $j(M)$ is linked with $\partial (h(B^k_r)\times B^n_R)$ with the linking number one. Consider the intersection
$V=X\times\R^n\cap N$. We may assume that $V$ is a manifold with boundary.
Since $\partial V$ is homologous to $\partial (h(B^k_r)\times B^n_R)$ in
$(\R^n\times\R^n)\setminus N$, the linking number of $\partial V$ and $j(M)$ is
one.
Since the intersection number of $V$ and $j(M)$ is one, the restriction
$h\mid_V:(V,\partial V)\to (B^{n+k}_d,\partial B^{n+k}_d)$ of a short map
$h:N\to B^{n+k}_d$ from Lemma 2 has degree one.\qed
\enddemo

\

\demo{Proof of Theorem 2} Let $\dim X=k$. Let $K$ be a triangulation on $X$
with diameters of simplices $\le 1$. We change the original embedding of $X$
into $l_2$ to the piecewise linear which is the same for vertices. We present
$X$ as a union of finite subcomplexes: $X=\cup K_i$, $K_i\subset K_{i+1}$.
Then every complex $K_i$ lies in a finite dimensional euclidean space
$\R^n_i\subset l_2$. Let $d$ be given. As in the prove of Theorem 1 we
construct a submanifold $V\subset X\times\R^{n(d)}$ and a short map of degree
one $f:(V,\partial V)\to (B^{k+n(d)}_d,\partial B^{k+n(d)}_d)$.
By Lemma 1 for large enough $r$ there is a retraction $\alpha$ of the $2d$-neighborhood
$N_{2d}(h(S_r^{k-1}))$ in $l_2$ onto $h(S_r^{k-1})$. There is $i$ such that
$h(S^{k-1}_r)\subset K_i$. Then we can work in $\R^{n_i}$ as in the
prove of Theorem 1 and we get $n(d)=n_i$. \qed
\enddemo

REMARK. If we replace in the above argument the Hilbert space $l_2$ by the Banach space $l_{\infty}$ we will obtain the following condition on $X$:
For every $\epsilon>0$ there exist $m$ and a submanifold with boundary
$W\subset X\times\R^m$ which admits an $\epsilon$-contracting map onto the
$l_{\infty}$ unit ball
$f:(W,\partial W)\to(B^{n+m}_{\infty},\partial B^{n+m}_{\infty})$ of degree
one. It is unclear if it would be possible to get a stable hypersphericity of $X$ from this (see [G4], page 8).  

\head \S3 Embedding Theorem for asymptotic dimension\endhead

We recall that the asymptotic dimension $asdim X$ of a metric space $X$ is
a minimal number $n$, if exists, such that for any $d>0$ there is a uniformly
bounded cover $\sU$ of $X$ which consists of $n+1$ $d$-disjoint families
$\sU=\sU^0\cup\dots\cup\sU^n$. A family of sets $\sV$ is $d$-disjoint if 
$d(V,V')=\inf\{d(x,x')\mid x\in V,x'\in V'\}>d$.
By $N_r(A)$ we denote the $r$-neigborhood if $r>0$ and the set
$A\setminus N_{-r}(X\setminus A)$ if $r\le 0$. Let $mesh\sU$ denote an upper bound for diameters of elements of a cover $\sU$ and let $L(\sU)$ denote the
Lebesgue number of $\sU$.
\proclaim{Proposition 1}
If a metric space $X$ with a base point $x_0$ has $asdim(X)\le n$ then there is a sequence of
uniformly bounded open covers $\sU_k$ of $X$, each cover $\sU_k$ splits into
a collection of $n+1$  $d_k$-disjoint families $\sU_k=\sU^0_k\cup\dots\cup\sU^n_k$
such that
\roster
\item{} $L(\sU_k)>d_k$ and $N_{-d_k}(U)\ne\emptyset$ for all $U\in\sU_k$,
\item{} $d_k>2^km_{k-1}$ where $m_{k-1}=mesh(\sU_{k-1})$,
\item{} For any $m\in \N$ and every $i\in\{0,\dots,n\}$ there is $l$ and $U\in\sU^i_l$ such that $N_{-d_k}(U)\supset B_m(x_0)$,
\item{} For every $i\in\{0,\dots,n\}$ and $k<l$ for any $U\in\sU_k^i$ and
$V\in\sU_l^i$ with $U\not\subset V$ there is the inequality $d(U,V)\ge 4$.
\endroster
\endproclaim
\demo{Proof}
We construct it by induction on $k$. We start with a cover $\sU_0$ with
$d_0>2$ and enumerate the partition $\sU_0=\sU^0_0\cup\dots\cup\sU^n_0$ in such a way that $d(x_0,X\setminus U)>d_0$ for some $U\in\sU^0_k$. We formulate the
condition (3) in a concrete fashion:

$(3)'$ For $l=m(n+1)+i$ where $i=l$ {\it mod} $n+1$ there is $U\in\sU^i_l$
such that $N_{-d_k}(U)\supset B_m(x_0)$.

Assume that the family $\{\sU_k\}$ is constructed for all $k\le l$ such that
the conditions $(1),(2),(3)',(4)$ hold. We define $d_{l+1}=2^{l+2}m_l$ and
consider a uniformly bounded cover $\bar\sU_{l+1}$ with the Lebesgue number
$L(\bar\sU_{l+1})>2d_{l+1}$
and with splitting in $n+1$ $d_{l+1}$-disjoint families
$\bar\sU_{l+1}=\bar\sU^0_{l+1}\cup\dots\cup\bar\sU^n_{l+1}$. Then we may assume that for all elements
$U\in\bar\sU_{l+1}$ we have $N_{-2d_{l+1}}(U)\ne\emptyset$. We just can delete
all elements $U$ from the cover $\bar\sU_{l+1}$ which do not satisfy that
property, and since $L(\bar\sU_{l+1})>2d_{l+1}$,  still we will have a cover of $X$. We enumerate families $\bar\sU^0_{l+1}\cup\dots\cup\bar\sU^n_{l+1}$
in such a way that $d(x_0,X\setminus U)>2d_{l+1}$ for some $U\in\bar\sU^i_{l+1}$
for $i=l+1$ {\it mod} $n+1$.
For every $U\in\bar\sU_{l+1}^i$ we define $\tilde U=U\setminus\bigcup_{V\not\subset U;V\in\sU^i_k,k\le l}\overline{N_4(V)}$.
We define $\sU^i_{l+1}=\{\tilde U\mid U\in\bar\sU^i_{l+1}\}$. Next we check all the properties.

(1). Take a point $x\in X$, then there exists $U\in\bar\sU_{l+1}$ such that
$B_{2d_{l+1}}(x)\subset U$. Then $B_{d_{l+1}}(x)\subset U\setminus N_{d_{l+1}}(X\setminus U)\subset U\setminus N_{m_l+4}(X\setminus U)\subset\tilde U$. Note that $N_{-d_{l+1}}(\tilde U)=\tilde U\setminus N_{d_{l+1}}(X\setminus U)\supset U\setminus N_{d_{l+1}+m_l+4}(X\setminus U)\supset N_{2d_{l+1}}(U)=
N_{-2d_{l+1}}(U)\ne\emptyset$.

(2). This condition holds by the definition.

$(3)'$. This condition holds by the construction.

(4). Assume that $V\in\sU^i_k$, $k\le l$ and $V\not\subset\tilde U$.
If $V\not\subset U$, then $d(V,\tilde U)\ge 4$. Assume that $V\subset U$.
Since $V\not\subset\tilde U$, there exists $V'\not\subset U$ and $V'\in\sU^i_s$
for some $s\le l$ such that $V\cap N_4(V')\ne\emptyset$. Hence $d(V,V')<4$.
By the condition (4) and induction assumption we have either $V\subset V'$ or
$V'\subset V$. In the first case it follows that $d(V,\tilde U)\ge d(V',\tilde U)\ge 4$. The second case is impossible, since $V'\subset V\subset U$ contradicts to the fact that $V'\not\subset U$. \qed 

\enddemo
\proclaim{Theorem 3}
Asume that $X$ is a metric space of bounded geometry with $asdim(X)\le n$. Then $X$ can be uniformly
embedded in a coarse sense in the product of $n+1$ locally finite trees.
\endproclaim
\demo{Proof}
Let $\sU_k$ be a sequence of covers of $X$ from Proposition 1. Let
$\sV_i=\cup_{k}\sU^i_k$. We define a map $\psi:\sV_i\to\sV_i$ by the following rule: $\psi(U)$ is the smallest $V\in\sV_i$ with respect to the inclusion
such that $V\ne U$ and $U\subset V$. The conditions 3) and 4) of Proposition 1 and
disjointness of $\sU^i_k$ for all $k$ imply that $\psi$ is well-defined.
For every $i\in\{0,\dots,n\}$ construct an oriented graph $T^i$ as follows.
For every $U\in\sV_i$ we consider an interval $I_U$ isometric with
$[0,2^k]$ and oriented from $2^k$ to $0$. For every $V\in\psi^{-1}(U)$ we
attach $I_V$ by the $0$-end to an integer point $a_V
=\min\{2^k,[\sup \phi_U(V)\frac{2^k}{d_k}]\}$ of $I_U$ where
$\phi_U(x)=d(x,X\setminus U)$ and $[a]$ means the integer part of $a$.
 
We show that the graph $T^i$ is a locally finite tree. For every $U,V\in\sV_i$
by the property 3) there exists $W\in\sV_i$ such that $U\cup V\subset W$.
This implies the connectedness of $T^i$. Since the orientation on $T^i$
defines a flow, i.e. every vertex is an initial point only for one arrow,
it follows that every cycle in $T^i$ must be oriented. Oriented cycles in
$T^i$ do not exist due to the inclusion nature of the orientation. Thus, $T^i$
is a tree. Note that for every nonzero vertex in $I_U$ only finitely many intervals are attached to it. It implies that if a vertex $v$ in $T^i$ is of infinite order, then
$v$ must be the $0$-vertex for all intervals involved. So a vertex $v$ of infinite order defines an infinite sequence $U_1\subset U_2\subset\dots\subset U_m\subset$ with $\psi(U_j)=U_{j+1}$ and $a_{U_j}=0$ for all $j$. Let $U_j\in\sU^i_{k_j}$, then
$\frac{2^{k_j}}{d_{k_j}}d(x,X\setminus U_{j+1})<1$ for all $x\in U_j$. Hence
$d(x,X\setminus U_{j+1}) < d_{k_j}$ for all $x\in U_j$, i.e.
$U_j\cap N_{-d_{k_j}}(U_{i+1})=\emptyset$. By the property 3) from Proposition 1
there is $U\in\sU^i_l$ with $N_{-d_l}(U)\supset U_1$. Then the condition 4)
implies that $U=U_j$ for some $j$ whence $l=k_j$. Therefore
$U_1\cap N_{-d_l}(U)=\emptyset$. Contradiction.

Next we define a map $p_i:X\to T^i$. By the condition 3) every point $x\in X$
is covered by some element $U\in\sU^i_k$ for some $k$. Let $U$ containing $x$ be taken
with the smallest $k$. We define $p_i(x)\in I_U$ as follows. Consider a map

\

$\xi:\bar N_{-d_k}(U)\cup\partial U\cup\bigcup_{V\in\psi^{-1}(U)}\partial V\to I_U=[0,2^k]$ 

\

defined as $\xi(\bar N_{-d_k}(U))=2^k$, 
 $\xi(\partial U)=0$ and $\xi(\partial V)=a_v$. Show that $\xi$ is a short map.
Let $y\in\partial V$ and $z\in\partial V'$,\  \ $V,V'\in\psi^{-1}(U)$ and
$V\ne V'$. Then by the condition 4) $d(y,z)\ge 4$. Note that

\

$|a_V-a_{V'}|\le\frac{2^k}{d_k}|d(y,X\setminus U)+m_{k-1}-d(z,X\setminus U)|+2
\le\frac{2^k}{d_k}d(y,z)+\frac{2^k}{d_k}m_{k-1}+2\le \frac{1}{4}d(y,z)+3\le
d(y,z)$. 

\

Here we applied the condition 2). Now let $y\in\partial V$ and
$z\in\partial U$. Then 

\

$|\xi(y)-\xi(z)|=a_V\le 
\frac{2^k}{d_k}(d(y,X\setminus U)+m_{k-1})\le\frac{1}{2m_{k-1}}d(y,z)+\frac{1}{2}\le d(y,z)$ 

\

provided
$d(y,z)>1$. Otherwise $a_v<1$ and hence $a_v=0$.
The case when $y\in\partial U$ and $z\in\bar N_{-d_k}(U)$ is obvious.
If $y\in\partial V$ and $z\in\bar N_{-d_k}(U)$ then

\

$|\xi(z)-\xi(y)|\le 2^k-\frac{2^k}{d_k}(d(y,X\setminus U)- m_{k-1}) +1\le
\frac{2^k}{d_k}(d_k-d(y,X\setminus U))+2\le \frac{2^k}{d_k}(d(z,X\setminus U)-
d(y, X\setminus U))+2\le d(y,z)$.

\

There exists a short extension $\bar\xi_U\to I_U$ of the map $\xi$.
We define $p_i(x)=\bar\xi_U(x)$. It easy to see that the map $p_i$ is short.

Show that the diagonal product $p=\Delta p_i:X\to\prod T^i$ is a uniform 
embedding. We consider $l_1$-metric on $\prod T^i$. Since each $p_i$ is short the map $p$ is Lipschitz. Clearly,  $Dist(p(x),p(x'))\ge \rho(d(x,x'))$ for
the function  $\rho(t)=\inf \{Dist(p(x),p(x'))\mid d(x,x')\ge t\}$.
 Assume that $\rho$ is bounded from above.
Then there is a sequence of pairs $(x_k,x'_k)$ of points with $d(x_k,x'_k)>m_k$ and
with $Dist(p(x_k),p(x'_k))\le b$ for all $k$.
For any $k$ there is an element $U\in\sU_k^i$ such that $d(x_k,X\setminus U)> d_k$.
 Since $d(x_k,x'_k)>m_k$, it follows that $x'_k\not\in U$. Note that
$\bar\xi_U(x_k)=2^k$ and hence the distance between $p_i(x_k)$ and
$p_i(x'_k)$ in the graph $T^i$ is greater than $2^k$. Therefore,
$Dist(p(x_k),p(x'_k))\ge 2^k$. Therefore $\rho$ tends to infinity and $p$ is
a uniform embedding.
\qed
\enddemo
\proclaim{Corollary}
Every metric space $X$ with $asdim(X)\le n$ is coarsely isomorphic to
a space $Y$ of a linear type.
\endproclaim

First we recall that according to Higson an asymptotically finite dimensional
space $Y$ has a {\it linear type} if in the above definition of $asdim$ 
there exists a constant $C$  such that for any $d$ the number $Cd$ is
an upper bound on the size of the cover $\sU$.
\demo{Proof}
Note that a tree and a finite product of trees are of linear type.
Hence every subspace of a finite product of trees is of linear type.
By Theorem 3 every asymptotically finite dimensional space $X$ is
coarsely isomorphic to a subset of a finite product of trees.
\enddemo

\proclaim{Lemma 3} Every locally finite tree is uniformly embeddable in
a complete simply connected 2-dimensional manifold $K(M)$ with negative curvature:
 $-k_1\le K(M)\le -k_2 < 0$.
\endproclaim
\demo{Proof}
Embed the tree into a plane. For every vertex $v$ which is not the end point we consider the longest
right-hand rule and left-hand rule paths from $v$.
If one of the paths is finite, then we end up in an end point of the tree. We
 attach a hyperbolic half-plane $H_v$ by isometry between the union of these paths and an interval 
in $\partial H_v$. For every end point vertex $v$ we will get two 
half-planes $H_{v_1}$ and $H_{v_2}$ with corresponding intervals ended in $v$. 
We attach them by isometries of
corresponding rays in $\partial H_{v_i}$. As the result we will get
a plane with a piecewise hyperbolic metric on it with possible 
singularities only at vertices of the tree. We can approximate this
metric by a smooth metric of strictly negative curvature
\enddemo
\proclaim{Theorem 4}
Every metric space $X$ with $asdim(X)\le n$ is uniformly embeddable in
a $(2n+2)$-dimensional non-positively curved manifold $W$ with 
$asdim(W)=2n+2$.
\endproclaim
\demo{Proof}
We apply Theorem 3 and Lemma 3 to obtain an embedding
of $X$ into $W=\prod M_i$ where each $M_i$ is 2-dimensional negatively curved
manifold. Then $W$ is non-positively curved. 
By a theorem of Gromov [G1] the asymptotic dimension of $M_i$ equals 2.
Hence the asymptotic dimension of $W$ is $2(n+1)$ [D-J] \qed
\enddemo
\proclaim{Problem 3}
Can $2n+2$ in the above theorem be improved to $2n+1$?
\endproclaim
Perhpas it is natural to ask whether every asymptotically $n$-dimensional
metric space of bounded geometry is embeddable into $(\H^2)^{n+1}$, the product of $n+1$ copies of the hyperbolic plane.

\head \S4 On the First Theorem of Yu\endhead

The following lemma is a generalization of Lemma 2.
\proclaim{Lemma 4}
Let $M$ be a closed $l$-dimensional manifold smoothly embedded in a
non-positively curved $n$-manifold $W$ with a trivial tubular $\epsilon$-neighborhood $N_{\epsilon}(M)$. Then for any $d>0$ there is
an embedding $\gamma:N_{\epsilon}(M)\to\R^n$ such that the diagonal
embedding $j=(1\Delta\gamma):M\to W\times\R^n$ has a regular neighborhood
$N$ with the projection $\nu: N\to j(M)$ such that:
\roster
\item{} there is a short map $h:N\to B^{2n-l}_d$ such that the restriction 
$h\mid_{\nu^{-1}(x)}:\nu^{-1}(x)\to B^{2n-l}_d$ is a homeomorphism for every $x\in j(M)$;
\item{} $pr_1(N)\subset N_{2d}(M)$ where $pr_1:W\times\R^n\to W$ is the
projection onto the first factor.
\endroster
\endproclaim
This Lemma together with the above Embedding Theorem allows to prove
the following:
\proclaim{Theorem 5}
If a uniformly contractible manifold $X$ has a finite asymptotic dimension,
then there is $n$  such that $X\times\R^n$ is integrally hyperspherical.
\endproclaim 
 The proof is exactly the same as in Theorem 1.
\demo{Proof of Lemma 4}
Let $TW$ be the tangent bundle of a non-positively curved complete simply 
connected $n$-dimensional Riemannian manifold $W$. 
For every $x\in W$ there is the exponential map $e_x:\R^n\to W$ which 
takes a vector $v$ to a point
$y=e_x(v)$ on the geodesic ray in the direction of $v$ with $d_W(x,y)=\|v\|$.
Note that $e_x$ is a homeomorphism for every $x$. The visual sphere at infinity
together with the exponential maps define a trivialization of $TW$.
A tubular $\epsilon$-neighborhood of a smmoth submanifold $M^l\subset W^n$ is a neighborhood $N_{\epsilon}(M)$ with the projection $p:N_{\epsilon}(M)\to M$
such that $p^{-1}(x)=e_x(B^{n-l}_{\epsilon})$ where $B^{n-l}_{\epsilon}\subset
N_x(M^l)\subset T_xW$ is an euclidean $\epsilon$-ball lying in the normal direction.

Using the embedding $e^{-1}_{x_0}\mid_{N_{\epsilon}(M)}:N_{\epsilon}(M)\to\R^n$ we define
$\gamma: N_{\epsilon}(M)\to\R^n$ to be a $\mu$-expanding map where $\mu$
is a large number defined as follows. By the definition of a tubular 
neighborhhod there is a 
lift $\e^{-1}:N_{\epsilon}(M)\to M\times B_{\epsilon}^{n-l}\subset TW$. 
Let $q:e^{-1}(N_{\epsilon}(M))\to B^{n-l}_{\epsilon}$ be a trivialization. 
Let $\lambda$ be a Lipschitz constant for the map $qe^{-1}$. 
Note that $e^{-1}(y)=(x,e^{-1}_x(y))\in T_xW$ for
$y\in p^{-1}(x)$. Let $N_{2d}(M)$ be a closed $2d$-neighborhood of $M$. The
correspondence $x\mapsto e^{-1}_x\mid_{N_{2d}(M)}$ defines a map
$\Phi:M\to C(N_{2d}(M),\R^n)$, where  the functional space 
$ C(N_{2d}(M),\R^n)$ is supplied with the {\it sup} norm 
$\|f\|=\sup_{z\in N_{2d}(M)}\|f(z)\|$. Let $\tilde\lambda>1$ be a 
Lipschitz constant of $\Phi$.
We take $\mu=\lambda\tilde\lambda d/\epsilon$. We define

\ \ \ \ \ \ \ \ \ \ \ \ \ \ \ \ \ \ \ \ \ \ \ \ \ \ \ \
$N=\bigcup_{x\in M} \bar B_{2d}(x)\times\gamma(\bar B_{\epsilon}(x))$.

Here we use $\bar B_r$ to denote a ball of radius $r$ in $W$ and $B_r$ for
a ball in $\R^n$. We define short maps $h_1:N\to B^n_d$ and $h_2:N\to 
B^{n-l}_d$ such that the sum $(h_1+h_2):N\to B^n_d\times B^{n-l}_d$ satisfies
the condition (1). We define $h_1\mid_{B_{2d}(x)\times\{y\}}=
\frac{1}{2}e^{-1}_x\mid_{B_{2d}(x)}$ for every
 $x\in M$ and every $y\in\gamma(B_{\epsilon}(x))$. Thus, $h_1((u,\gamma(v))_x)=
\frac{1}{2}e^{-1}_x(u)$ for $u\in B_{2d}(x)$ and $v\in B_{\epsilon}(x)$.
We define $h_2=\frac{d}{\epsilon}q\circ e^{-1}\circ\gamma^{-1}\circ pr_2$.
Then we can estimate a Lipschitz constant for $h_2$ as the product
$\frac{d}{\epsilon}\lambda\mu^{-1}=\frac{1}{\tilde\lambda}\le 1$. 
Let $z=(u,\gamma(v))_x$ and $z'=(u',\gamma(v'))_{x'}$ be two points in $N$.
Then $\|h_1(z)-h_1(z')\|=\frac{1}{2}\|e^{-1}_x(u)-e^{-1}_{x'}(u')\|\le
\frac{1}{2}\|e^{-1}_x(u)-e^{-1}_{x'}(u)\|+\frac{1}{2}\|e^{-1}_{x'}(u)-
e^{-1}_{x'}(u')\|\le\frac{1}{2}\tilde\lambda d(x,x')+\frac{1}{2}d(u,u')\le
\frac{1}{2}(\mu d(x,x')+d(u,u'))\le\frac{\sqrt{2}}{2}
\sqrt{d^2(u,u')+d^2(\gamma(v),\gamma(v'))}
\le d(z,z')$.
Then we define $h$ as the composition of $h_1+h_2$ and the natural projection of
$B^n_d\times B^{n-l}_d$ onto $B^{2n-l}_d$. The condition (1) holds.
The condition (2) holds by the definition of $N$. \qed
\enddemo

\enddocument